\documentclass[twocolumn]{autart}    

\usepackage{graphicx}          
\usepackage{color}
\usepackage{graphicx}
\usepackage{amsfonts,amssymb}
\usepackage{bm}
\usepackage{bbm}
\usepackage{amsmath}
\usepackage{mathrsfs}
\usepackage{subcaption}

\newcommand{\blkdiag}{\text{blkdiag}}
\newcommand{\ibLambda}{{\it{\bar \Lambda }}}
\newcommand{\ihLambda}{{\it{ \hat \Lambda}}}
\newcommand{\blue}{\color{blue}}

\definecolor{grey}{rgb}{0.55,0.55,0.55}

\begin{document}

\begin{frontmatter}

\title{Balanced Truncation of Networked Linear Passive Systems
	\thanksref{footnoteinfo}} 

\thanks[footnoteinfo]{Corresponding author: Xiaodong~Cheng.}

\author[Cheng1,Cheng2]{Xiaodong Cheng}\ead{x.cheng@rug.nl},    
\author[Cheng1]{Jacquelien M.A. Scherpen}\ead{j.m.a.scherpen@rug.nl},
\author[Cheng1]{Bart Besselink}\ead{b.besselink@rug.nl}  

\address[Cheng1]{Jan C.
	Willems Center for Systems and Control, Faculty of Science and Engineering, University
	of Groningen, The Netherlands.}  
\address[Cheng2]{Control Systems Group, Department of Electrical Engineering, Eindhoven University of Technology, The Netherlands}


\begin{abstract}                          

This paper studies model order reduction of multi-agent systems consisting of identical linear passive subsystems, where the interconnection topology is characterized by an undirected weighted graph. Balanced truncation based on a pair of specifically selected generalized Gramians is implemented on the asymptotically stable part of the full-order network model, which leads to a reduced-order system preserving the passivity of each subsystem. 
Moreover, it is proven that there exists a coordinate transformation to convert the resulting reduced-order model to a state-space model of Laplacian dynamics. Thus, the proposed method simultaneously reduces the complexity of the network structure and individual agent dynamics, and it preserves the passivity of the subsystems and the synchronization of the network. Moreover, it allows for the \textit{a priori} computation of a bound on the approximation error. Finally, the feasibility of the method is demonstrated by an example.

\end{abstract}

\end{frontmatter}

\section{Introduction}

Multi-agent systems, or network systems, recently have become a rapidly evolving area of research with a tremendous amount of applications, including power grids, cooperative robots, biology and chemical reaction networks (see, e.g., \cite{mesbahi2010graph,ren2005survey} for an overview). However, a multi-agent system may be high-dimensional due to the large scale of networks and complexity of nodal dynamics. In most cases, the full-order complex network models are neither practical nor necessary for controller design, system simulation and validation. Hence, it is desirable to apply model order reduction techniques to derive a lower-order
approximation of the original network system with an acceptable accuracy. 

In many network applications, Laplacian structures play an important role, as they represent communication graphs characterizing the interactions among agents. For instance, the synchronization and stability of networks are analyzed in the context of Laplacian dynamics (see, e.g., \cite{consensus2010li,mesbahi2010graph,ren2005survey,Monshizadeh2013stability}). Thus, it is a natural requirement to preserve the algebraic structure of the Laplacian matrix in order to inherit a network interpretation in a reduced-order model, where a reduced Laplacian matrix is employed to describe diffusive coupling protocols of the reduced network.

Conventional reduction techniques, including balanced truncation, Hankel-norm approximation, and Krylov subspace methods, do not explicitly take the interconnection structure into account in deriving the reduced-order models. Consequently, the direct application of these methods to multi-agent systems potentially leads to the loss of desired properties such as the synchronization of networks and the structure of the subsystems.  Towards the model reduction with the preservation of network structure, mainstream methodologies are focusing on graph clustering. From, e.g., \cite{Schaft2014,Monshizadeh2014,ishizaki2015clustereddirected,XiaodongTAC2018MAS,XiaodongTAC20172OROM,XiaodongACOM2018Power}, we have observed that the clustering-based approaches naturally maintain the spatial structure of networks and show an insightful interpretation for the reduction process. Nevertheless, the approximation of these methods relies on the selection of clusters, while finding a reduced network with the smallest error generally is an NP-hard problem, see \cite{SurveyClustering}. For tree networks, \cite{Besselink} considers the so-called edge dynamics of a network, where a pair
of diagonal generalized Gramian matrices of the edge system are used for characterizing the importance of the edges. Then, the nodes linked by less important edges are clustered, resulting in an \textit{a priori} bound on the approximation error. However, the application of this approach is only applicable to a tree topology. Another method based on singular perturbation is developed for reducing network complexity, which is mainly applied to electrical grids and chemical reaction networks (see e.g., \cite{Chow2013PowerReduction,Monshizadeh2018constant,Rao2013graph}). In these works, the network
structure is preserved as the Schur complement of the Laplacian matrix of the original network is again a Laplacian matrix, representing a
smaller-scale network. This approach is of particular interest for simplifying networked single integrators, while it may be less suitable for dealing with networks of higher-order agent dynamics since the Laplacian is not the only matrix any more defining the network dynamics.
The other direction in model order reduction of multi-agent systems is to reduce the dimension of each individual subsystem, see e.g., \cite{Monshizadeh2013stability,sandberg2009interconnected,XiaodongECC2018Lure}, which use the generalized balanced truncation method to reduce subsystems in a network, while keeping the interconnection topology untouched.  

In this paper, we {
develop a technique that can reduce the complexity of network structures and individual agent dynamics \textit{simultaneously}, extending preliminary results in \cite{XiaodongIFAC2017}. This problem setting has seldomly been studied in the literature so far, and different from \cite{Ishizaki2016dissipative}, we aim to reduce the network structure and agent dynamics in a unified framework.
Particularly, this paper considers multi-agent systems composed of identical higher-order linear passive subsystems, where the network topology is characterized by an undirected weighted graph. 
The core step in the proposed technique is balancing the asymptotically stable part of the network system based on generalized Gramians. 
After truncating the balanced model, we obtain a reduced-order system with a lower dimension, which preserves the passivity of the subsystems. Although the network structure is not necessarily preserved in this step, 
we show that there exists a set of coordinates in which the reduced-order model can be interpreted as a network system. Specifically, the main contributions of this paper are summarized as follows. First, two generalized Gramians that are structured by the Kronecker product are selected such that the balanced truncation is applied to reduce the network structure and agent dynamics in a unified framework. The proposed method guarantees an \textit{a priori} computation of a bound on the approximation error with respect to external inputs and outputs. Second, we propose a necessary and sufficient condition of a matrix being similar to a Laplacian matrix (see Theorem~\ref{thm:LapReal}). With this result, the reduction process is designed to preserve the Laplacian structure in the reduced network.

The remainder of this paper is organized as follows. Section \ref{sec:ProblemFormulation} provides necessary preliminaries and formulates the
model reduction problem of networked passive systems. Then, Section \ref{sec:MainResults} presents the main results, that is the model reduction procedure for a network system. The proposed method is illustrated by an example in Section \ref{sec:Example}, and finally conclusions are made in Section \ref{sec:Conclusion}.

\textit{Notation:} The symbol $\mathbb{R}$  denotes the set of real numbers, whereas $I_n$ and $\bm{1}_n$ represent the identity matrix of size $n$ and all-ones vector of $n$ entries, respectively. The 2-induce norm of matrix $A$ is denoted by $\| A \|_2$.
The Kronecker product of matrices $A \in \mathbb{R}^{m \times n}$ and $B\in \mathbb{R}^{p \times q}$ is denoted by $A \otimes B \in \mathbb{R}^{mp \times nq}$. 
Besides, $\bm{\Sigma}$ represents a linear system, and the operation $\bm{\Sigma}_1 +\bm{\Sigma}_2$ means the parallel interconnection of two linear systems by summing their transfer functions. The $\mathcal{H}_\infty$-norm of the transfer function of a linear system $\bm{\Sigma}$ is denoted by $\lVert \bm{\Sigma} \rVert_{\mathcal{H}_\infty}$.


\section{Preliminaries and Problem Formulation}
\label{sec:ProblemFormulation}
%
 
Consider a network of $N$ nodes, and the dynamics on each node is described by 
\begin{equation} \label{sysagent}
{\bm{\Sigma}_i}: \ \left \{
\begin{array}{l}
\dot{x}_i = A x_i + B \nu_i, \\
\eta_i = C x_i,
\end{array}
\right.	
\end{equation}
where $x_i \in \mathbb{R}^{n}$, $\nu_i\in \mathbb{R}^{m}$ and $\eta_i \in \mathbb{R}^{m}$ are the states, control inputs and outputs of agent $i$, respectively. Throughout the paper, we assume that the system realization in (\ref{sysagent}) is \textit{minimal} and \textit{passive}. 
Passivity is a natural property of many real physical systems, including mechanical systems, power networks, and thermodynamical systems (see \cite{Koeln2017GraphPassivity,Hatanaka2015Passivity}). The passivity of ${\bm{\Sigma}_i}$ can be charaterized by the following lemma.   
\begin{lem} \cite{Willems1972Dissipative,Schaft2008BalanPassive} \label{lem:KYP}
	The linear system ${\bm{\Sigma}_i}$ in (\ref{sysagent}) is passive if and only if there exists a symmetric positive definite matrix $K$ such that 
	\begin{equation} \label{eq:KYP}
		A^\top K + K A \leq 0,  \ C = B^\top K.
	\end{equation}
\end{lem}

The agents are assumed to interact with each other through a weighted undirected connected graph $\mathcal{G}$ containing $N$ nodes. More precisely, we have the following diffusive coupling rule 
\begin{equation} \label{eq:protocol}
	\nu_i = -\sum_{j=1,j\ne i}^{N} w_{ij} \left(\eta_i-\eta_j\right)
+ \sum_{j=1}^{p} f_{ij} u_j,
\end{equation}
where $w_{ij} \geq 0$ represents the strength of the coupling between nodes $i$ and $j$. Moreover, $u_j \in \mathbb{R}^{m}$ with $j=\{1,2, \cdots,p\}$ are external inputs, and $f_{ij} \in \mathbb{R}$ is the amplification of the $j$-th input acting on agent $i$, which is zero when $u_j$ has no effect on node $i$. 
Let $ 
	y_i = \sum_{j=1}^{N} h_{ij} \eta_j, \ i=1,2,\cdots,q,
$
with $h_{ij} \in \mathbb{R}$, be the the $i$-th external output. We then obtain the overall multi-agent system in a compact form:
\begin{equation} \label{sysh}
\bm{\Sigma}:  
\left \{
\begin{array}{l}
\dot{x} = \left(I_N \otimes A - L \otimes BC\right)x + (F \otimes B)u,
\\
y=(H \otimes C)x.
\end{array}
\right.
\end{equation}
Here, $F\in \mathbb{R}^{N \times p}$ and $H \in \mathbb{R}^{q \times N}$ are the collections of $f_{ij}$ and $h_{ij}$, respectively, and
$
	x := \left[x_1^\top, \cdots, x_n^\top \right]^\top \in \mathbb{R}^{N n}, 
	u := \left[{u}_1^\top, \cdots, u_p^\top \right]^\top \in \mathbb{R}^{p m},
	y := \left[{y}_1^\top, \cdots, {y}_q^\top \right]^\top \in \mathbb{R}^{q m}
$.
Furthermore, $L \in \mathbb{R}^{N \times N}$ is the \textit{Laplacian matrix} of the graph $\mathcal{G}$ with the $(i,j)$-th entry as
\begin{equation} \label{defn:Laplacian}
L_{ij} = \left\{ \begin{array}{ll} 
\sum_{j=1,j\ne i}^{N} w_{i j}, & \text{if} \ i = j,\\
-w_{i,j}, & \text{otherwise.}
\end{array}
\right.
\end{equation}
In this paper,  the underlying graph $\mathcal{G}$ is assumed to be \textit{undirected} (i.e., $L_{ij} = L_{ji}$) and \textit{connected}, in which case $L$ has the following properties, see, e.g., \cite{XiaodongTAC20172OROM,XiaodongTAC2018MAS}.
\begin{rem} \label{rem:StruCond}	
	For a connected undirected graph, the Laplacian matrix $L$ fulfills the following \textit{structural conditions}:
	(1) $\bm{1}_N^\top L = 0$, and $L \bm{1}_N = 0$;
	(2) $L_{ij} \leq 0$ if $i \ne j$, and $L_{ii} > 0$;
	(3) ${L}$ is positive semi-definite with a single zero eigenvalue.	
	 The Laplacian $L$ is the matrix representation of the graph $\mathcal{G}$. Conversely, a real square matrix can be interpreted as a connected undirected graph if it satisfies the above conditions.
\end{rem}
Laplacian matrices are commonly used for describing network systems with diffusive couplings and are very instrumental for the 
synchronization analysis of networks. For the network in (\ref{sysh}), the synchronization property is characterized in the following lemma.
\begin{lem}   \label{lem:syn} \cite{Besselink}
	Consider the network system $\bm{\Sigma}$ in (\ref{sysh}). If the graph $\mathcal{G}$ is connected, and each subsystem ${\bm{\Sigma}_i}$ in (\ref{sysagent}) is observable, then $\bm{\Sigma}$  synchronizes for $u = 0$, i.e.,
	\begin{equation}
	\lim\limits_{t \rightarrow \infty}  \left[x_i(t)-x_j(t)\right]   = 0, \ \ \forall i,j \in \{1,2,\cdots,N\}.
	\end{equation}
	for any initial condition $x_i(0)$, $i = 1, 2, \cdots, N$.
\end{lem}

%

{
}

Now, we address the model order reduction problem for multi-agent systems of the form (\ref{sysh}) as follows.

\begin{prob} \label{prob:appx}
	Given a multi-agent system $\bm{\Sigma}$ as in (\ref{sysh}), find a reduced-order model
	\begin{equation} \label{sysrh0}
	\bm{\hat{\Sigma}}: \left \{
	\begin{array}{l}
	\dot{\hat{x}} = (I_{k} \otimes \hat{A} - \hat{L} \otimes \hat{B}\hat{C})\hat{x} + (\hat{F} \otimes \hat{B}) u, \\
	\hat{y} = (\hat{H} \otimes \hat{C}) \hat{x},
	\end{array}
	\right.
	\end{equation}
	such that the following objectives are achieved:
	
	\begin{itemize}
		\item $\hat{L} \in \mathbb{R}^{k \times k}$, with $k \leq N$, is an undirected graph Laplacian satisfying the structural conditions in Remark \ref{rem:StruCond}.  
		
		\item The lower-order approximation of the agent dynamics
		\begin{equation} \label{sysagentRed}
		\bm{\hat{\Sigma}}_i: \left \{
		\begin{array}{l}
		\dot{\hat{x}}_i = \hat{A} \hat{x}_i + \hat{B} \hat{\nu}_i, \\
		\hat{\eta}_i = \hat{C} \hat{x}_i,
		\end{array}
		\right.	
		\end{equation}
		with the reduced state vector $\hat{x}_i \in \mathbb{R}^{r}$ ($r \leq n$), is passive, i.e., satisfies the condition in Lemma~\ref{lem:KYP}. 
		
		\item The overall approximation error $\lVert \bm{\Sigma}-\bm{\hat{\Sigma}} \rVert_{\mathcal{H}_\infty}$ is small.  
	\end{itemize}
%
%
\end{prob} 


\section{Main Results}
\label{sec:MainResults}


\subsection{Separation of Network System}
\label{sec:Separation}

Since $A$ in (\ref{sysagent}) is not necessarily Hurwitz, meaning that  $\bm{\Sigma}$ may be not asymptotically stable,  a direct application of   balanced truncation   to $\bm{\Sigma}$ is not feasible. We thereby introduce a decomposition of $\bm{\Sigma}$ using 
 the following \textit{spectral decomposition} of the graph Laplacian as  
\begin{equation} \label{eq:LEigDec}
L = T {\it \Lambda} T^\top = \begin{bmatrix}
T_1 & T_2
\end{bmatrix} 
\begin{bmatrix}
\ibLambda & \\  & 0
\end{bmatrix}
\begin{bmatrix}
T_1^\top \\ T_2^\top
\end{bmatrix},
\end{equation}
where $T_2 = \bm{1}_N/\sqrt{N}$ and
$ \label{eq:barLambda}
	\ibLambda : = \text{diag}(\lambda_1, \lambda_2, \cdots, \lambda_{N-1}),
$
with $\lambda_1 \geq \lambda_2 \geq \cdots \geq \lambda_{N-1} > 0 $ the nonzero eigenvalues of $L$.
Then, {\blue we
apply the coordinate transformation $x = ({T} \otimes I) z$ to the system $\bm{\Sigma}$}, which yields two independent components, namely, an \textit{average module}
\begin{equation} \label{sysa}
\bf{\Sigma_a}: 
\left\{
\begin{array}{lr}
\dot{z}_a = Az_a + \dfrac{1}{\sqrt{N}}  (\bm{1}_N^\top F\otimes B) u, \\
y_a = \dfrac{1}{\sqrt{N}}  ( H \bm{1}_N \otimes C) z_a, \\	
\end{array}
\right.
\end{equation}
with ${\blue z_a : = (\bm{1}^\top_N/\sqrt{N} \otimes I) x} \in \mathbb{R}^{n}$,
and an asymptotically stable system
\begin{equation} \label{syst}
\bf{\Sigma_s}: 
\left\{
\begin{array}{lr}
\dot{z}_s = (I_{N-1} \otimes A - \ibLambda \otimes BC) z_s + ( \bar{F} \otimes B) u, \\
y_s = (\bar{H} \otimes C) z_s.\\	
\end{array}
\right.
\end{equation}
where ${\blue z_s: = (T_1^\top \otimes I) x} \in \mathbb{R}^{(N-1)\times n}$, $\bar{F}=T_1^{T} F$, {and} $\bar{H} = H T_1$. Note that the synchronization property of $\bm{\Sigma}$ implies  the asymptotic stability of $\bf{\Sigma_s}$, see  \cite{Besselink}.
Thus, we can apply
balanced truncation to $\bf{{\Sigma}_s}$ to generate its lower-order approximation $\bf{\hat{\Sigma}_s}$. It meanwhile gives a reduced subsystem  $(\hat{A}, \hat{B}, \hat{C})$ resulting in a reduced-order average module $\bf{\hat{\Sigma}_a}$. Combining  $\bf{\hat{\Sigma}_s}$ with $\bf{\hat{\Sigma}_a}$ formulates a reduced-order model $\bm{\tilde{\Sigma}}$ whose input-output behavior is similar to that of the original system $\bm{\Sigma}$. 
However, at this stage, the network structure is not necessarily preserved by $\bm{\tilde{\Sigma}}$. Then, it is desired to use a coordinate transformation to convert $\bm{\tilde{\Sigma}}$ to $\hat{\bm{\Sigma}}$, which restores the Laplacian structure. The whole procedure is summarized in Fig. \ref{procedure}, and the detailed implementations are discussed in the following subsections.

\begin{figure}[!tp]\centering
	\centering
	\includegraphics[width=0.4\textwidth]{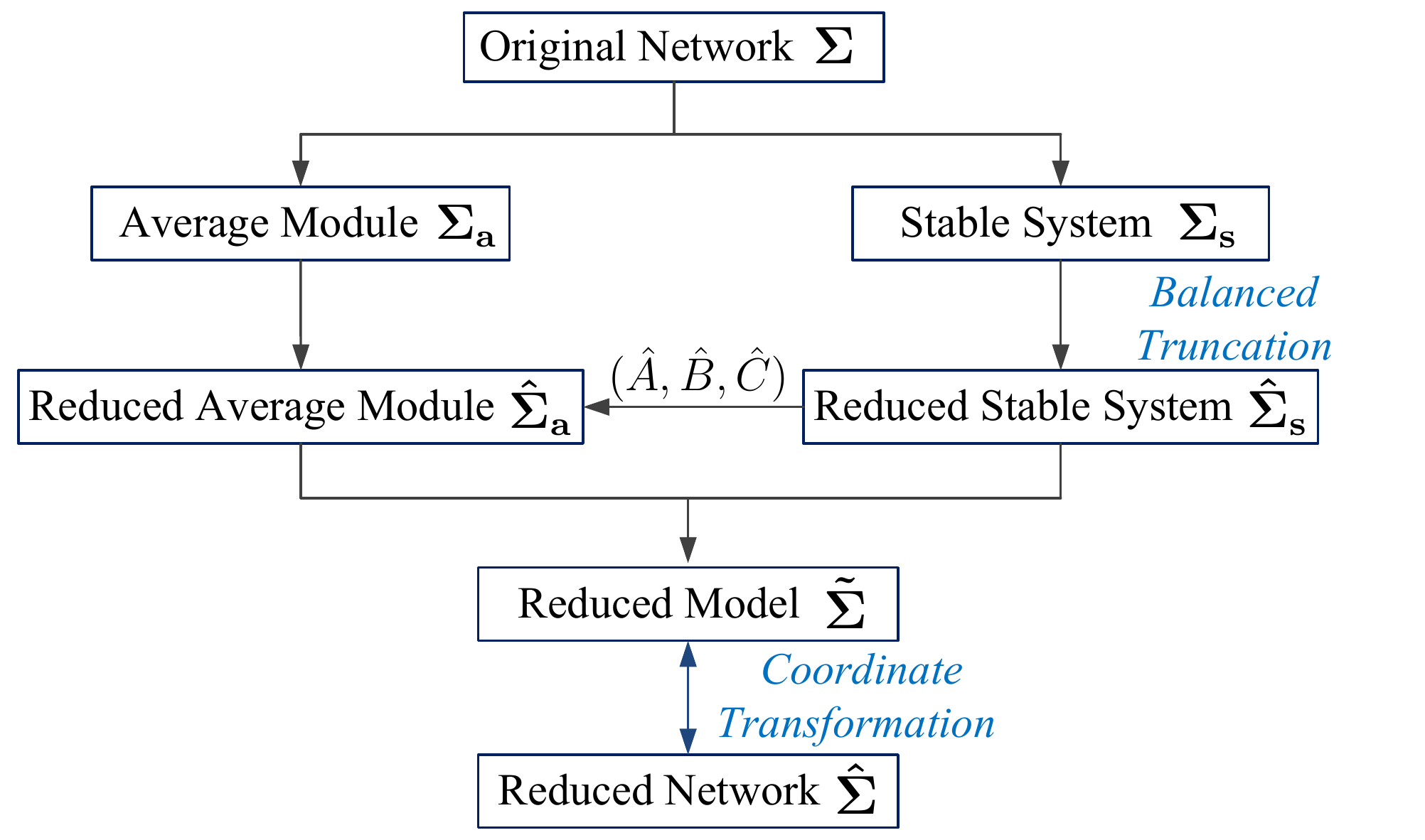}	
	\caption{The scheme for the structure preserving model order reduction of networked passive systems}	
	\label{procedure}
\end{figure}

\subsection{Balanced Truncation by Generalized Gramians}
\label{sec:Balancing}

{ 
Following \cite{Dullerud2013CourseRobust}, the generalized Gramians of the asymptotically stable system $\bf{\Sigma_s}$ are defined. 
\begin{defn} \label{defn:GenGrams}
	Consider the stable system $\bf{\Sigma_s}$ and denote $\Phi := I \otimes A - \ibLambda \otimes BC$. Two positive definite matrices $\mathcal{X}$ and $\mathcal{Y}$ are said to be the {generalized controllability and observability Gramians} of \ $\bf{\Sigma_s}$, respectively, if they satisfy 
	\begin{subequations}
		\begin{align}
		\Phi \mathcal{X} + \mathcal{X} \Phi^\top + (\bar{F} \otimes B) (\bar{F}^\top \otimes B^\top) &\leq {0}, \label{eq:ConGramInq}
		\\ 
		\Phi^\top \mathcal{Y} + \mathcal{Y} \Phi + (\bar{H}^\top \otimes C^\top) (\bar{H} \otimes C) &\leq {0}. \label{eq:ObsGramInq}
		\end{align}
	\end{subequations}
	Moreover, a {generalized balanced realization} is achieved when $\mathcal{X} = \mathcal{Y} > 0$ are diagonal. The diagonal entries are called  generalized Hankel singular values  (GHSVs).	
\end{defn}
}


Suppose $\ibLambda$ in (\ref{eq:LEigDec}) has $s$ distinct diagonal entries ordered as: $\bar{\lambda}_1 > \bar{\lambda}_2 > \cdots > \bar{\lambda}_s$. We rewrite $\ibLambda$ as
$ 
\ibLambda = \blkdiag (\bar{\lambda}_1 I_{m_1}, \bar{\lambda}_2 I_{m_2} \cdots, \bar{\lambda}_s I_{m_s}),
$
where $m_i$ is the multiplicity of $\bar{\lambda}_i$, and $\sum_{i=1}^{s} m_i = N-1$. 
Then, we consider  
the following Lyapunov equation and inequality:
\begin{subequations} \label{eq:Lyap}
	\begin{align}
	-\ibLambda {X}  - {X}  \ibLambda + \bar{F}\bar{F}^\top &= 0, \label{eq:LyapX}
	\\ 
	-\ibLambda {Y}  - {Y} \ibLambda + \bar{H}^\top\bar{H} &\leq 0, \label{eq:LyapY} 
	\end{align}
\end{subequations}
where ${X} = {X}^\top > 0$ and  
$  
	{Y} := \blkdiag (Y_1, Y_2, \cdots, Y_s),
$
with $Y_i=Y_i^\top> 0$ and $Y_i \in \mathbb{R}^{m_i \times m_i}$, for $i = 1,2, \cdots, s$. 
The block-diagonal structure of $Y$ is crucial to guarantee  that the reduced-order model, obtained by preforming balanced truncation on the basis of $X$ and $Y$, to be interpreted as a network system again, see Lemma~\ref{lem:PosiEigen} and Theorem~\ref{thm:LapReal}. The matrix $X$ is chosen as the standard controllability Gramian for a smaller error bound. 
  Compared with our former notation in \cite{XiaodongIFAC2017}, the Definition~of the observability Gramian is more general, since it is not necessary to be strictly diagonal.

\begin{rem} \label{rem:XY}
	There exist a variety of networks, especially symmetric ones such as stars, circles, chains or complete graphs, whose Laplacian matrices have repeated eigenvalues.   
	Particularly, when $L$ refers to a complete graph with identical weights, all the eigenvalues in are equal, leading to a full matrix $Y$, and (\ref{eq:LyapY}) can be specialized to an equality. Besides, by the duality between controllability and observability, we can also use $-\ibLambda {X}  - {X}  \ibLambda + \bar{F}\bar{F}^\top \leq  0$, and $-\ibLambda {Y}  - {Y} \ibLambda + \bar{H}^\top\bar{H} = 0$ to characterize the pair $X$ and $Y$ for the balanced truncation, where now $X$ is constrained to have a block-diagonal structure.
\end{rem}

{The existence of the solutions $X$ and $Y$  in (\ref{eq:LyapX}) and (\ref{eq:LyapY}) are guaranteed, as $\ibLambda>0$ is positive diagonal. Furthermore, in practice, the generalized observability Gramian is obtained by minimizing the trace of $Y$, see, e.g., \cite{Besselink,sandberg2009interconnected}.} Based on $X$ and $Y$, we further define a pair of generalized Gramians for the stable system $\bf{\Sigma_s}$.

%
%
%

\begin{thm} \label{thm:GenGram}
	Consider $X$, $Y$ as the solutions of (\ref{eq:Lyap}), and let $K_{m}>0$ and $K_{M}>0$ be the minimum and maximum solutions of (\ref{eq:KYP}). 
	Then, the matrices
	\begin{equation} \label{eq:GenGrams}
	\mathcal{X} := {X} \otimes K_{M}^{-1} \ \text{and} \ \mathcal{Y}: =  {Y} \otimes K_m
	\end{equation} 
	characterize generalized Gramians of the asymptotically stable system $\bf{\Sigma_s}$.
	Moreover, there exist two nonsingular matrices $T_\mathcal{G}$ and $T_\mathcal{D}$ such that $\mathcal{T} = T_\mathcal{G} \otimes T_\mathcal{D}$ satisfies
	\begin{equation} \label{eq:balancing}
	\mathcal{T} \mathcal{X} \mathcal{T}^\top = \mathcal{T}^{-T}\mathcal{Y}\mathcal{T}^{-1} = \Sigma_\mathcal{G} \otimes \Sigma_\mathcal{D}.
	\end{equation}
	Here, $\Sigma_\mathcal{G} := \mathrm{diag} \{\sigma_1, \sigma_2,\cdots, \sigma_{N-1}\},$ and $\Sigma_\mathcal{D} := \mathrm{diag} \{\tau_1, \tau_2,\cdots, \tau_{n}\}$, where $\sigma_1 \geq \sigma_2 \geq \cdots \geq \sigma_{N-1},$ and $\tau_1 \geq \tau_2 \geq \cdots \geq \tau_{n}
	$ are
	corresponding to the square roots of the spectrum of $XY$ and $K_{M}^{-1}  K_{m}$, respectively. 
\end{thm}
\begin{pf}
    By the passivity of  ${\bm{\Sigma}_i}$ and Lemma~\ref{lem:KYP},
   	we verify that 
	\begin{equation*}
	\begin{split}
	&\Phi \mathcal{X} + \mathcal{X} \Phi^\top + (\bar{F} \otimes B) (\bar{F}^\top \otimes B^\top) \\
	=& X \otimes (A K_{M}^{-1} + K_{M}^{-1} A^\top) \\
	 & \hspace{1cm}+ (-\ibLambda {X}  - {X}  \ibLambda + \bar{F}\bar{F}^\top) \otimes BB^\top \leq 0, \\
	\end{split}
	\end{equation*}
	where the inequality holds due to (\ref{eq:LyapX}) and $K_{M}$ being a solution of (\ref{eq:KYP}). Similarly, it can be verify that $\mathcal{Y}$ in (\ref{eq:GenGrams}) satisfies the inequality in (\ref{eq:ObsGramInq}).
	Thus, by Definition~\ref{defn:GenGrams}, $\mathcal{X}$ and $\mathcal{Y}$ in (\ref{eq:GenGrams}) characterize the generalized Gramians of $\bf{\Sigma_s}$. 
	Next, by the standard balancing theory \cite{antoulas2005approximation}, there exist nonsingular matrices $T_\mathcal{G}$ and $T_\mathcal{D}$ such that 
	\begin{subequations}  
	\begin{align}
	T_\mathcal{G} X T_\mathcal{G}^\top 
	 &=  \Sigma_{\mathcal{G}} = T_\mathcal{G}^{-T} Y T_\mathcal{G}^{-1}, \label{eq:TG}
	 \\
		T_\mathcal{D} K_{M}^{-1} T_\mathcal{D}^\top &= \Sigma_{\mathcal{D}} = T_\mathcal{D}^{-T} K_{m} T_\mathcal{D}^{-1}. \label{eq:TD}
	\end{align} 
	\end{subequations}
	Thus, $ 
	\mathcal{T} = T_\mathcal{G} \otimes T_\mathcal{D},
	$ can be used for the balancing transformation of $\bf{\Sigma_s}$. Moreover, since
	$T_\mathcal{G} X Y T_\mathcal{G}^{-1} = 
	T_\mathcal{G} X T_\mathcal{G}^\top T_\mathcal{G}^{-T} Y T_\mathcal{G}^{-1} = 
	\Sigma_{\mathcal{G}}^2$ and $T_\mathcal{D} K_{M}^{-1} K_{m} T_\mathcal{D}^{-1} = 
	T_\mathcal{D} K_{M}^{-1} T_\mathcal{D}^\top T_\mathcal{D}^{-T} K_{m} T_\mathcal{D}^{-1} 
	= \Sigma_{\mathcal{D}}^2$, 
	the singular values in $\Sigma_\mathcal{G}$ and $\Sigma_\mathcal{D}$ are characterized by the square roots of the spectrum of $XY$ and $K_{M}^{-1}  K_{m}$, respectively. 
\end{pf}

\begin{rem}
		The maximum and minimum solutions of (\ref{eq:KYP}), $K_M$ and $K_m$, respectively define 
		 the \textit{available storage} $
		\frac{1}{2} \langle x, K_{m}x \rangle
		$ and the \textit{required supply}
		$ \frac{1}{2} \langle x, K_{M}x \rangle $
		of the agent system \cite{Willems1972Dissipative}.
		Any $K>0$ satisfying (\ref{eq:KYP}) will lie between these two extremal
		solutions, i.e., $0 < K_m \leq K \leq K_M$.
		It is also noted that the solution of (\ref{eq:KYP}) may be unique, i.e., $K_M = K_m$, e.g., when the system (\ref{sysagent}) is \textit{lossless} \cite{Schaft2008BalanPassive} or $B$ is square and nonsingular \cite{Willems1972Dissipative}. In this case, we have $\Sigma_{\mathcal{D}} = I_n$ meaning that the subsystems are not suitable for reduction. If $K_{M}  \ne K_{m}$, it can be verified that the diagonal entries of $\Sigma_{\mathcal{D}}$ in (\ref{eq:balancing}) satisfy $\tau_i \leq 1$, $\forall i=1,2,\cdots,n$.
\end{rem}

Generally, there exist multiple choices of generalized Gramians as the solutions of (\ref{eq:ConGramInq}) and (\ref{eq:ObsGramInq}). This paper specifically selects the pair of Gramians in (\ref{eq:GenGrams}) with the Kronecker product structure, implying that they can be simultaneously diagonalized, (i.e., balanced) using transformations of the form $\mathcal{T} = T_\mathcal{G} \otimes T_\mathcal{D}$. Note that $T_\mathcal{G}$ and $T_\mathcal{D}$   are independently generated from (\ref{eq:Lyap}) and (\ref{eq:KYP}). More precisely, $T_\mathcal{G}$ only depends on the network structure, or the triplet $(\ibLambda, \bar{F}, \bar{H})$, while $T_\mathcal{D}$ only replies on the  agent dynamics, i.e., the triplet $(A, B, C)$. Thus, the Laplacian dynamics and the subsystem (\ref{sysagent}) can be reduced independently, allowing the resulting reduced-order model to preserve a network interpretation as well as the passivity of subsystems. Denote
\begin{equation}
	(\ihLambda_{1}, \hat{F}_{1}, \hat{H}_{1})\ \text{and}\  (\hat{A}, \hat{B}, \hat{C}): = {\bm{\hat{\Sigma}}_i} 
\end{equation}
as the reduced-order models of $(\ibLambda, \bar{F}, \bar{H})$ and $(A, B, C)$, respectively, where 
$\ihLambda_{1} \in \mathbb{R}^{(k-1) \times (k-1)}$,
$\hat{F}_{1} \in \mathbb{R}^{(k-1) \times  p}$, 
$\hat{H}_{1} \in \mathbb{R}^{q \times (k-1)}$, $\hat{A} \in \mathbb{R}^{r \times r}$,
$\hat{B} \in \mathbb{R}^{r \times  m}$, and
$\hat{C} \in \mathbb{R}^{m \times r}$.  
Consequently, the reduced-order models of the average module (\ref{sysa}) and the stable system (\ref{syst}) are constructed.
\begin{subequations}  
\begin{align} 
\label{sysaRed}
\bf{\hat{\Sigma}_a}: &
\left\{
\begin{array}{lr}
\dot{\hat{z}}_a = \hat{A} \hat{z}_a + \dfrac{1}{\sqrt{N}}  (\bm{1}_N^\top F\otimes \hat{B}) u, \\
\hat{y}_a = \dfrac{1}{\sqrt{N}}  ( H \bm{1}_N \otimes \hat{C}) \hat{z}_a. \\	
\end{array}
\right.
\\
\label{systRed}
\bf{\hat{\Sigma}_s}: &
\left\{
\begin{array}{lr}
\dot{\hat{z}}_s = (I_{k-1} \otimes \hat{A} - \ihLambda_{1} \otimes \hat{B}\hat{C}) \hat{z}_s  + (\hat{F}_{1} \otimes \hat{B}) u, \\
\hat{y}_s = 
(\hat{H}_{1} \otimes \hat{C}) \hat{z}_s.\\	
\end{array}
\right.
\end{align}
\end{subequations}  
\begin{rem}
	When ${\bm{\Sigma}_i}$ is strictly passive \cite{Opmeer2013GapMetric}, the balanced truncation of ${\bm{\Sigma}_i}$ on the basis of {\blue $K_m$ and $K_M$} delivers a strictly passive and minimal reduced-order model $\bm{\hat{\Sigma}}_i$. Otherwise, $\bm{\hat{\Sigma}}_i$ is passive but not necessarily minimal. Nevertheless, we can always replace $\bm{\hat{\Sigma}}_i$ by its minimal realization $(\hat{A}, \hat{B}, \hat{C})$ as in \cite{Polyuga2008Kalman}, and the replacement does not change the transfer functions of $\bf{\hat{\Sigma}_s}$ and $\bf{\hat{\Sigma}_a}$. 
\end{rem}

Combining the reduced-order models  $\bf{\hat{\Sigma}_a}$ and $\bf{\hat{\Sigma}_s}$ formulates a lower-dimensional approximation of  $\bm{\Sigma}$ as 
\begin{equation} \label{sysRedMod}
\bm{\tilde{\Sigma}}: 
\left\{
\begin{array}{lr}
\dot{\hat{z}} = (I_k \otimes \hat{A} - \mathcal{N} \otimes \hat{B}\hat{C}) \hat{z}  + (\mathcal{F} \otimes \hat{B}) u, \\
\hat{y} = 
(\mathcal{H} \otimes \hat{C}) \hat{z}.\\	
\end{array}
\right.
\end{equation}
where 
$$
\mathcal{N} =  \begin{bmatrix}
\ihLambda_{1} & \\ & 0
\end{bmatrix}, 
\mathcal{F} = \begin{bmatrix}
\hat{F}_{1} \\  \tfrac{1}{\sqrt{N}}  \bm{1}_N^\top F 
\end{bmatrix},
\mathcal{H} = \begin{bmatrix}
\hat{H}_{1} & \tfrac{1}{\sqrt{N}} H \bm{1}_N 
\end{bmatrix}.
$$
Here, $\mathcal{N}$ is not yet a Laplacian matrix, which prohibits the interpretation of $\bm{\tilde{\Sigma}}$ as a network system. However, $\mathcal{N}$ has the following property.
\begin{lem} \label{lem:PosiEigen}
	The matrix $\mathcal{N}$ in (\ref{sysRedMod}) has only one zero eigenvalue at the origin and all the other eigenvalues are positive and real.	
\end{lem}
\begin{pf}
	Using the structure property of $Y$, we verify that $Y \ibLambda = Y^{1/2} \ibLambda Y^{1/2}$.
	The reduced matrix $\ihLambda_{1}$ in (\ref{sysRedMod}) is obtained by the following standard projection 
	\begin{align} \label{eq:VQVVQ}
		\nonumber \ihLambda_{1} &= \left[(V_1^\top Y V_1)^{-1}  V_1^\top Y\right]  \ibLambda V_1 \\
		&= (V_1^\top Y V_1)^{-1}  V_1^\top Y^{1/2} \ibLambda Y^{1/2} V_1,
	\end{align}
	where $V_1 \in \mathbb{R}^{N \times k}$ is the left projection matrix obtained by the singular value decomposition of $X^{1/2} Y^{1/2}$, see \cite{antoulas2005approximation} for more details. As $V_1$ is full column rank,  (\ref{eq:VQVVQ}) shows that $\ihLambda_{1}$ is the product of two positive definite matrices, implying that $\ihLambda_{1}$ only has positive and real eigenvalues.
\end{pf}
\begin{rem} 
	Generally, balanced truncation does not preserve the realness of eigenvalues. Lemma~\ref{lem:PosiEigen} is the result of using a generalized observability Gramian $Y$ with the block diagonal structure. As mentioned in Remark~\ref{rem:XY}, we may exchange the equality and inequality in (\ref{eq:Lyap}) because of
	duality. Then, the eigenvalue realness of $\ihLambda_{1}$ is also guaranteed due to the similar reasoning.
\end{rem}}

\subsection{Network Realization}
\label{sec:Realization}

The spectral property of $\mathcal{N}$ allows for a reinterpretation of the reduced-order model $\bm{\tilde{\Sigma}}$  as a network system again. 
\begin{thm} \label{thm:LapReal}
	A real square matrix $\mathcal{N}$ is similar to a Laplacian matrix $\mathcal{L}$ associated with an undirected connected graph, if and only if $\mathcal{N}$ is diagonalizable and has exactly one zero eigenvalue while all the other eigenvalues are real positive.
\end{thm}
The proof is provided in the appendix.
By Theorem~\ref{thm:LapReal}, we can achieve a network realization of $\bm{\tilde{\Sigma}}$, and at least a complete network is guaranteed to be obtained. 
Specifically, we find 
a nonsingular matrix $\mathcal{T}_n$ such that
$
\hat{L} = \mathcal{T}_n^{-1} \mathcal{N} \mathcal{T}_n,
$
where $\hat{L}$ is Laplacian matrix characterizing a reduced connected undirected graph with $k$ nodes.
Applying the coordinate transform $\hat{z} = (\mathcal{T}_n \otimes I_{{r}}) \hat{x}$ to  $\bm{\tilde{\Sigma}}$ in (\ref{sysRedMod}) then yields a reduced-order network model
\begin{equation} \label{sysRedNet}
\bm{\hat{\Sigma}}: 
\left\{
\begin{array}{lr}
\dot{\hat{x}} = (I_k \otimes \hat{A} - \hat{L} \otimes \hat{B}\hat{C}) \hat{x}  + (\hat{F} \otimes \hat{B}) u, \\
\hat{y} = 
(\hat{H} \otimes \hat{C}) \hat{x},\\	
\end{array}
\right.
\end{equation}
with
$
\hat{F} = \mathcal{T}_n^{-1} \mathcal{F}  \ \text{and} \
\hat{H} = \mathcal{H} \mathcal{T}_n.
$
Since the reduced subsystem $(\hat{A},\hat{B},\hat{C})$ is passive and minimal, the following theorem holds.
\begin{thm}
	The reduced networked passive system $\bm{\hat{\Sigma}}$ in (\ref{sysRedNet}) preserves synchronization, i.e., when $u=0$, it holds that
	\begin{equation}
	\lim\limits_{t \rightarrow \infty}  \left[\hat{x}_i(t)-\hat{x}_j(t)\right]  = 0, \ \ \forall i,j \in \{1,2,\cdots,k\},
	\end{equation}
	for any initial condition $\hat{x}(0)$.
\end{thm}

\subsection{Error Analysis}

Following the separation of the multi-agent system $\bm{\Sigma}$ in Section \ref{sec:Separation}, we analyze the approximation error for the overall system by using the triangular inequality. 
\begin{align} \label{eq:ErrParts}
	\nonumber \lVert \bm{\Sigma} - \bm{\hat{\Sigma}} \rVert_{\mathcal{H}_\infty} 
	&= \lVert (\bf{\Sigma_s} + \bf{\Sigma_a}) - (\bf{\hat{\Sigma}_s}  + \bf{\hat{\Sigma}_a}) \rVert_{\mathcal{H}_\infty} \\
	&\leq 
	\lVert \bf{\Sigma_s} - \bf{\hat{\Sigma}_s} \rVert_{\mathcal{H}_\infty} +
	\lVert \bf{\Sigma_a} - \bf{\hat{\Sigma}_a} \rVert_{\mathcal{H}_\infty}.
\end{align}
First, an \textit{a priori} bound on the approximation error of the stable system $\bf{\Sigma_s}$ is provided. 
\begin{lem}\label{lem:errboundstab}
    Consider the stable system $\bf{\Sigma_s}$ in (\ref{syst}) and its approximation $\bf{\hat{\Sigma}_s}$ in (\ref{sysaRed}). The approximation error has an upper bound as 
	$
		\lVert \bf{\Sigma_s} - \bf{\hat{\Sigma}_s} \rVert_{\mathcal{H}_\infty}  \leq \gamma,
	$
	where 
	\begin{equation}\label{eq:errboundstab}
		\gamma = 2 \sum\limits_{i=k}^{N-1}\sum\limits_{j=1}^{n}\sigma_i \tau_j + 2\sum\limits_{i=1}^{k-1}\sum\limits_{j=r+1}^{n}\sigma_i \tau_j.
	\end{equation}
	with $\sigma_i$ and $\tau_i$ the diagonal entries of $\Sigma_{\mathcal{G}}$ and  $\Sigma_{\mathcal{D}}$ in (\ref{eq:balancing}), respectively.
\end{lem}
\begin{pf}
	The GHSVs of the balanced system of $\bf{\Sigma_s}$ are ordered on the diagonal of $\Sigma_\mathcal{G} \otimes \Sigma_\mathcal{D}$ as  \begin{equation*} 
		\begin{split}
		&\Sigma_\mathcal{G} \otimes \Sigma_\mathcal{D} =\\ & \blkdiag \left(	\sigma_1 \begin{bmatrix}
		\tau_1 &  & \\
		& \ddots & \\
		&  & \tau_n\\
		\end{bmatrix}, \cdots, \sigma_{N-1} \begin{bmatrix}
		\tau_1 &  & \\
		& \ddots & \\
		&  & \tau_n\\
		\end{bmatrix} \right).
		\end{split}
		\end{equation*} 
	Then, the bound $\gamma$ is obtained from the standard error analysis for balanced truncation. 
\end{pf}

The approximation error on the average module, i.e., $\bf{\Sigma_a}-\bf{\hat{\Sigma}_a}$,  is given by
\begin{equation*} \label{eq:Delta}
\begin{split}
\Delta_a(s) = &\dfrac{1}{N} (H\bm{1}_N \otimes C) (sI_n-A)^{-1} (\bm{1}_N^\top F \otimes B) \\
& - \dfrac{1}{N} (H\bm{1}_N \otimes \hat{C}) (sI_r-\hat{A})^{-1} (\bm{1}_N^\top F \otimes \hat{B}) 
\\
  = &\dfrac{H\bm{1}_N \bm{1}_N^\top F}{N} \otimes \Delta_i(s),
\end{split}
\end{equation*}
where 
$\Delta_i(s)$ 
is the transfer function of ${\bm{\Sigma}_i}-\bm{\hat{\Sigma}}_i$. Hence, the approximation error on the average module is bounded if and only if the error between the original and reduced agent dynamics is bounded. Note that $\bm{\hat{\Sigma}}_i$ is obtained from positive real balancing of ${\bm{\Sigma}_i}$, and generally, there does not exist an \textit{a priori} bound on $\lVert {\bm{\Sigma}_i} - \bm{\hat{\Sigma}}_i \rVert_{\mathcal{H}_\infty}$. Nevertheless, \textit{a posteriori} bound can be obtained, see \cite{Opmeer2013GapMetric}. If  $\Delta_i(s) \in \mathcal{H}_\infty$, we obtain 
\begin{equation} \label{eq:errboundaverage}
\begin{split}
&\lVert {\bf{\Sigma_a}} - {\bf{\hat{\Sigma}_a}} \rVert_{\mathcal{H}_\infty} \leq \dfrac{\gamma_a}{N} \lVert {\bm{\Sigma}_i} - \bm{\hat{\Sigma}}_i \rVert_{\mathcal{H}_\infty} 
\end{split}
\end{equation}
with $\gamma_a: = \lVert H\bm{1}_N \bm{1}_N^\top F \rVert_2$. 


In the rest of this section, special cases are discussed where an \textit{a priori} upper bound on $\lVert \bm{\Sigma} - \bm{\hat{\Sigma}} \rVert_{\mathcal{H}_\infty}$ in (\ref{eq:ErrParts}) can be obtained.
The first case is when we only reduce the dimension of the network while the agent dynamics are untouched as in \cite{Besselink}. In this case, we obtain $\lVert \mathbf{\Sigma_a} - \mathbf{\hat{\Sigma}_a} \rVert_{\mathcal{H}_\infty} = 0$, which yields the error bound straightforwardly following from Lemma~\ref{lem:errboundstab}. 
\begin{thm}
	Consider the network system $\bm{\Sigma}$ with $N$ agents and its reduced-order model $\bm{\hat{\Sigma}}$ with $k$ agents. If the agent system ${\bm{\Sigma}_i}$ is not reduced, the error bound
	\begin{equation}
		\lVert \bm{\Sigma} - \bm{\hat{\Sigma}} \rVert_{\mathcal{H}_\infty} = 
		\lVert {\bf{\Sigma_s}} - {\bf{\hat{\Sigma}_s}} \rVert_{\mathcal{H}_\infty} \leq 2 \sum_{i=k}^{N-1}\sum_{j=1}^{n} \sigma_i \tau_j,
	\end{equation}
	holds, where $\sigma_i$ and $\tau_i$ are defined in Theorem~\ref{thm:GenGram}.
\end{thm}

The second case is when  the average module is not observable from the outputs of the overall system $\bm{\Sigma}$ or uncontrollable by the external inputs. Specifically, we have 
\begin{equation} \label{eq:specialcase2}
	H\bm{1}_N = 0, \ \text{or} \ \bm{1}_N^\top F = 0,
\end{equation}
which also implies
$\lVert \mathbf{\Sigma_a} - \mathbf{\hat{\Sigma}_a} \rVert_{\mathcal{H}_\infty} = 0$.
In practice, this means that we only observe or control the differences between the agents. 
Such differences usually play a crucial role in distributed control of networks, which aims to steer the  states of (partial) nodes to achieve a certain agreement.
A typical example can be found in \cite{Monshizadeh2013stability,Monshizadeh2014} where $H$ in (\ref{sysh}) is the \textit{incidence matrix} of the underlying network. 
\begin{cor} \label{coro:ErrorBound}
	Consider the network system $\bm{\Sigma}$ with $N$ agents and its reduced-order network model $\bm{\hat{\Sigma}}$ with $k$ agents. If $H\bm{1}_N = 0$ or $\bm{1}_N^\top F = 0$, the approximation between $\bm{\Sigma}$ and $\bm{\hat{\Sigma}}$ is bounded by
	\begin{equation*}
	\lVert \bm{\Sigma} - \bm{\hat{\Sigma}} \rVert_{\mathcal{H}_\infty}
	= 	\lVert \bf{\Sigma_s} - \bf{\hat{\Sigma}_s} \rVert_{\mathcal{H}_\infty}
	\leq  \gamma,
	\end{equation*}
	where
	$\gamma$ is given in (\ref{eq:errboundstab}).
\end{cor}

\section{Illustrative Example}
\label{sec:Example}

\begin{figure}[!tp]\centering
	\begin{minipage}[t]{0.55\linewidth}
		\centering
		\includegraphics[width=0.7\textwidth]{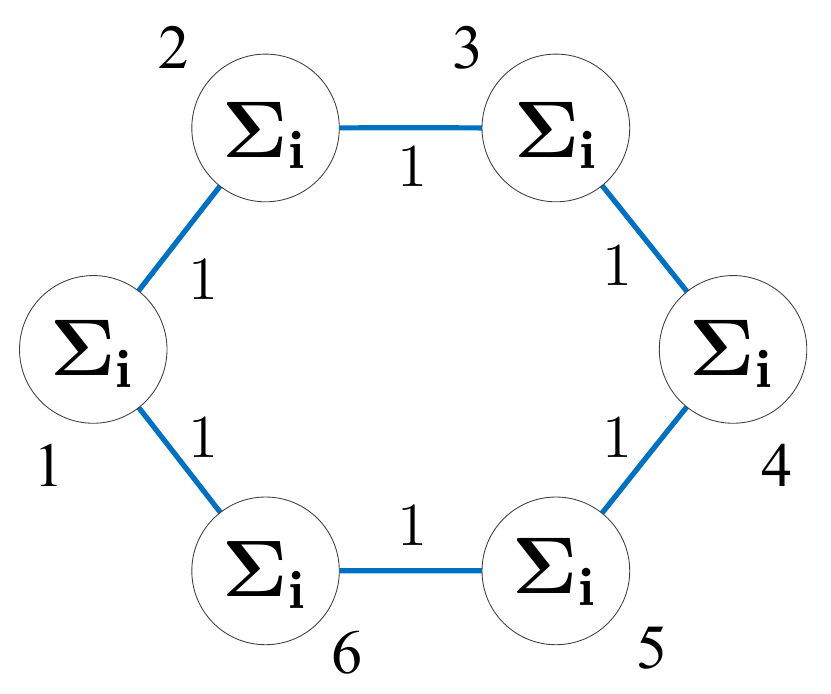}
		\subcaption{}
		\label{fig:originalnetwork}		
	\end{minipage}%
	\hfill
	\begin{minipage}[t]{0.42\linewidth}
		\centering
		\includegraphics[width=0.75\textwidth]{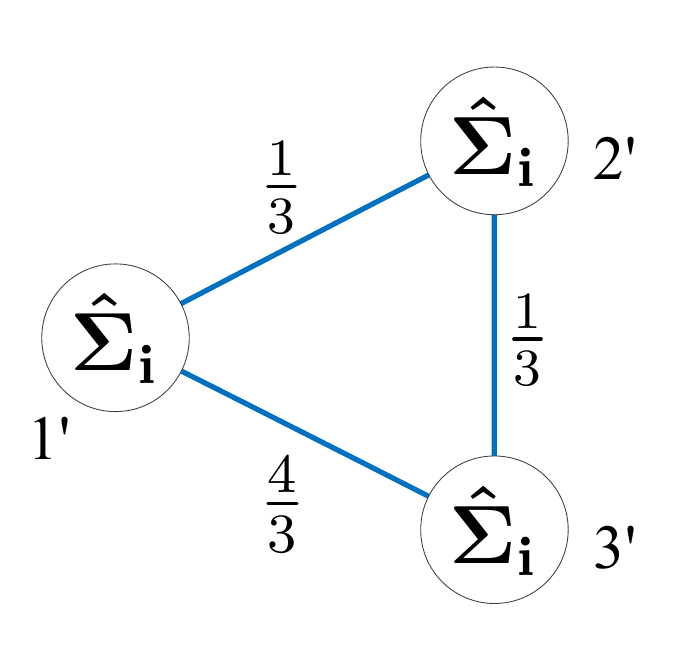}	
 		\subcaption{}
 		\label{fig:reducednetwork}
	\end{minipage}
\caption{(a) and (b) illustrate the original and reduced communication graph, respectively.}
\end{figure}

To demonstrate the feasibility of the proposed method, we consider networked robotic manipulators as a multi-agent system example.  The dynamics of each rigid robot manipulator is described as a standard mechanical system in the form (\ref{sysagent}) with 
\begin{equation}\label{sys:portHam}
	A = \begin{bmatrix}
	\bm{0} & M^{-1}\\-I & -DM^{-1}
	\end{bmatrix}, \
	C = B^\top \cdot \begin{bmatrix}
	I & \bm{0} \\ \bm{0} & M^{-1}
	\end{bmatrix},
\end{equation}
where $D \geq 0$ and $M>0$ are the system damping and mass-inertia
matrices, respectively. By Lemma~\ref{lem:KYP}, each manipulator agent is passive since there exists a positive definite matrix $P: = \blkdiag(I, M^{-1})$ satisfying (\ref{eq:KYP}). In this example, the system parameters in (\ref{sys:portHam}) are specified as $M = \frac{1}{2} I_4$, $B =[
0, 0, 0, 0, 1, 0, 0, 0 ]^\top$, and
{\scriptsize
\begin{equation*}
	D = \begin{bmatrix}
	2  &  -1  &   0  &   0\\
	-1 &    4  &  -2  &   0\\
	0  &  -2  &   4  &  -1\\
	0 &    0  &  -1  &   2\\
	\end{bmatrix}.
\end{equation*}
}%
which yields the dynamics of each individual agent  with state dimension $n=8$. 
Furthermore, $6$ agents communicate according to an undirected cyclic graph depicted in Fig. \ref{fig:originalnetwork}.
The Laplacian matrix and external input and output matrices are given by
{\scriptsize
\begin{equation*}
\begin{split}
	 &L = \begin{bmatrix}
	 2  &  -1  &  0  &   0   &  0  &  -1\\
	 -1  &   2  &  -1  &   0  &   0  &   0 \\
	 0  &  -1  &   2 &   -1   &  0   &  0\\
	 0 &    0  &  -1  &   2  &  -1  &   0\\
	 0  &   0  &   0  &  -1  &  2  &  -1\\
	 -1   &  0  &   0  &   0  &  -1 &   2
	 \end{bmatrix}, \
	 F = \begin{bmatrix}
	 1 \\0.5\\0\\0\\0\\0
	 \end{bmatrix}, \ 
	 H^\top = \begin{bmatrix}
	 1 \\ 0 \\ -1 \\ 0 \\ 0 \\ 0
	 \end{bmatrix}.
\end{split}
\end{equation*}
}%
It can be verified that the subsystems ${\bm{\Sigma}_i}$ is minimal. Thus, the overall network is synchronized  by Lemma~\ref{lem:syn}.

The nonzero eigenvalues of $L$ are $\lambda_1 = 4, \lambda_2 = \lambda_3 = 3, \lambda_4 = \lambda_5 = 1$. Solving the linear matrix inequality (\ref{eq:LyapY}) by minimizing the trace of $Y$, we obtain  
{\scriptsize
	\begin{equation*}
	Y = \blkdiag \left( 
	\begin{bmatrix}
	0.0120  &  0.0964 \\
	0.0964 &   0.7766
	\end{bmatrix},
	\begin{bmatrix}
	0.3416  &  0.1972   \\
	0.1972  &  0.1139
	\end{bmatrix},
	2.23 \cdot 10^{-5}
	\right).
	\end{equation*}}%
Moreover, from (\ref{eq:LyapX}) and (\ref{eq:KYP}), we compute $X$, $K_M$ and $K_m$, respectively. In this example, $K_M \ne K_m$ holds.

\begin{figure}[!tp]\centering
	\centering
	\includegraphics[width=0.5\textwidth]{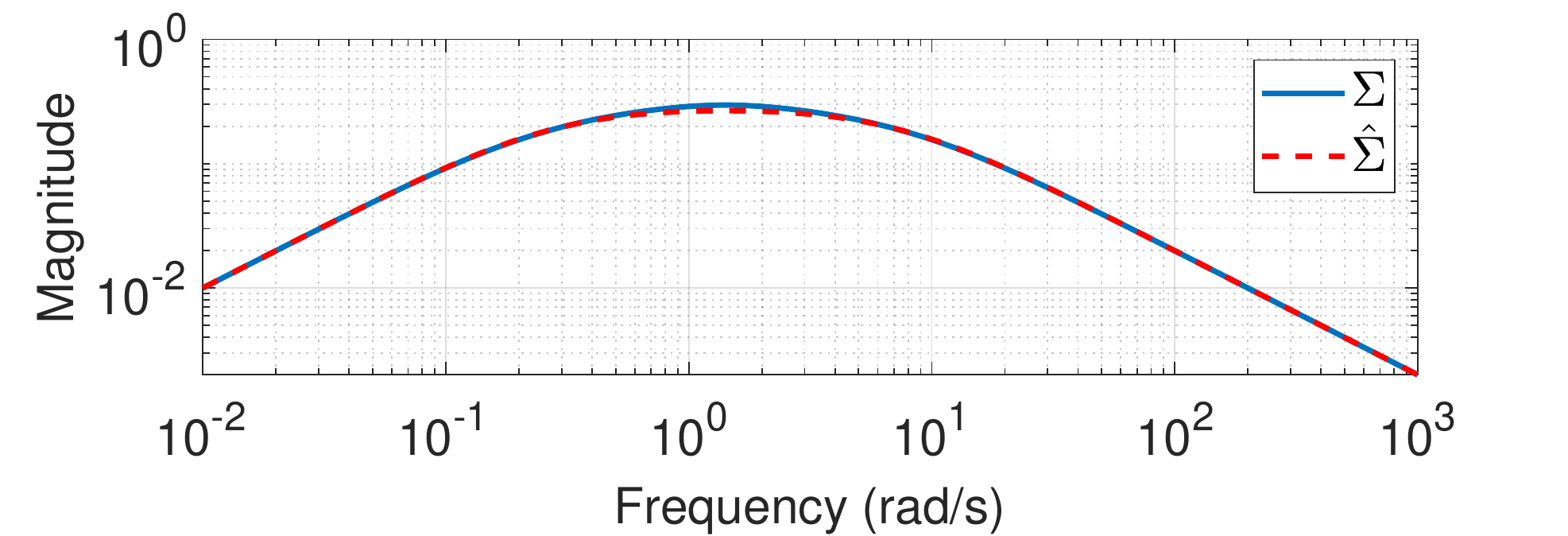}	
	\caption{The frequency responses of the original and reduced multi-agent systems, which are represented by the solid and dashed lines in the plot respectively.}	
	\label{bode}
\end{figure}

The goal is to reduce the dimension of the agent systems to $r = 2$ and the number of nodes to $k=3$.
Applying the generalized balanced truncation discussed in Section~\ref{sec:Balancing}, we obtain a reduced-order subsystem $\bm{\hat{\Sigma}}_i$ with
{\scriptsize
\begin{equation*}
\hat{A} = \begin{bmatrix}
   0  & -1.4142\\
	1.4142  & -4
\end{bmatrix}, \ 
\hat{B} = \hat{C}^\top = \begin{bmatrix}
   0\\
	-1.4142
\end{bmatrix}.
\end{equation*}
}%
Furthermore, by the network realization method in Section \ref{sec:Realization}, a lower-dimensional Laplacian matrix and external input and output matrices can be computed as
{\scriptsize
	\begin{equation*}
	\begin{split}
	\hat{L} &= \dfrac{1}{3} \begin{bmatrix}
	 5 &  -1 &  -4 \\
	-1  &   2 &   -1\\
	-4  & -1 &    5
	\end{bmatrix}, \ 
	\hat{F} = \begin{bmatrix}
	 -0.9270\\
	1.1380\\
	0.8496
	\end{bmatrix}, \ 
	\hat{H}^\top  = \begin{bmatrix}
	-0.4939 \\   0.4249 \\   0.0690
	\end{bmatrix}.
		\end{split}
	\end{equation*}
}%
Note that $\hat{L}$ represents a reduced interconnection network as shown in Fig. \ref{fig:reducednetwork}, which consists of $3$ reduced agents. We observe that $\bm{\hat{\Sigma}}_i$ is passive and minimal. Therefore, the reduced-order multi-agent system preserves the synchronization property.
Next, to compare the input-output behavior of the reduced-order network to the original one, we plot the frequency responses of both systems in Fig. \ref{bode} and compute the actual model reduction error: $\lVert \bm{\Sigma} - \bm{\hat{\Sigma}} \rVert_{\mathcal{H}_\infty} \approx 0.0295$. Since $H \bm{1}_6 = 0$, we then obtain the \textit{a priori} error bound by Corollary \ref{coro:ErrorBound} as $\lVert \bm{\Sigma} - \bm{\hat{\Sigma}} \rVert_{\mathcal{H}_\infty} \leq 0.0773$. 
Therefore, the original network is well approximated by the reduced-order model.

\section{Concluding Remarks}
\label{sec:Conclusion}

In this paper, we have developed a novel structure-preserving model reduction method for networked passive systems.  
Based on the selected generalized Gramians, the dimension of each subsystem and the network topology are reduced via a unified framework of balancing.
The resulting model is guaranteed to be converted to reduced-order network system.  Moreover, an \textit{a priori} error bound on the overall system has been provided. 	For future works, multi-agent systems with nonlinear agent dynamics and communication protocols are of interest.

\appendix  
\section{Proof of Theorem~\ref{thm:LapReal}}
\label{Append:LapReal}
\vspace{-20pt}
\begin{pf}
		The ``only if'' part can be seen from Remark \ref{rem:StruCond}. The rest of the proof shows the ``if'' part.
		Let $\mathcal{N} \in \mathbb{R}^{n \times n}$ be diagonalizable, and denote its eigenvalues as
		\begin{equation} \label{eq:lamdas}
		\lambda_1 \geq \lambda_2 \geq \cdots \geq \lambda_{n-1} > \lambda_{n}=0.
		\end{equation}
		Then, there exists a spectral decomposition $\mathcal{N} = T_1 D_1 T_1^{-1}$ with $D_1 = \mathrm{diag}(\lambda_1, \lambda_2, \cdots, \lambda_{n})$. 
		
		On the other hand, any undirected graph Laplacian $\mathcal{L}$ can be written in the form of 
		{\scriptsize
			\begin{equation} \label{eq:GeneralFormLaplacian}
			\mathcal{L} = \left[\begin{matrix}
			\alpha_1 & -w_{1,2} & \cdots & -w_{1,n} \\
			-w_{2,1} & \alpha_2  & \cdots & -w_{2,n} \\
			\vdots & \vdots & \ddots & \vdots \\
			-w_{n,1} & -w_{n,2} & \cdots & \alpha_n
			\end{matrix}\right],
			\end{equation}}%
		where $w_{i,j} = w_{j,i} \geq 0$ denotes the weight of edge $(i,j)$, which is the same as $w_{ij}$ in (\ref{eq:protocol}), and 
		\begin{equation} \label{eq:alpha}
		\alpha_i = \sum_{j=1,j \ne i}^{n} w_{i,j}.
		\end{equation}
		There exists an eigenvalue decomposition $\mathcal{L} = T_2 D_2 T_2^{-1}$. If $D_1 = D_2$,
		the following equation holds
		\begin{equation} \label{eq:SimNL}
		\mathcal{L} = (T_2 T_1^{-1}) \mathcal{N} (T_2 T_1^{-1})^{-1}.
		\end{equation}
		Hence, it is sufficient to prove that there always exists a set of weights $w_{i,j}$ such that the resulting Laplacian matrix $\mathcal{L}$ in (\ref{eq:GeneralFormLaplacian}) and $\mathcal{N}$ have the same eigenvalues (\ref{eq:lamdas}).

		Consider the characteristic polynomial of $\mathcal{L}$, i.e., 
		{\scriptsize
			\begin{equation*}
			\left|  \mathcal{L} -  \lambda I_n \right| = 
			\left|  \begin{matrix}
			\alpha_1-\lambda & -w_{1,2} & \cdots & -w_{1,n-1} & -w_{1,n} \\
			-w_{1,2} & \alpha_2-\lambda  & \cdots & -w_{2,n-1} & -w_{2,n} \\
			\vdots & \vdots & \ddots & \vdots & \vdots\\
			-w_{1,n-1} & -w_{2,n-1} & \cdots &\alpha_{n-1}-\lambda & -w_{n-1,n}\\
			-w_{1,n} & -w_{2,n} & \cdots & -w_{n-1,n} &\alpha_n-\lambda
			\end{matrix}\right|.
			\end{equation*}}%
		As elementary row operations do not change the determinant, we sum all rows to the final row to obtain
		{\scriptsize
			\begin{equation*}
			\left|  \mathcal{L} -  \lambda I_n \right| = 
			\left|  \begin{matrix}
			\alpha_1-\lambda & -w_{1,2} & \cdots & -w_{1,n-1} & -w_{1,n} \\
			-w_{1,2} & \alpha_2-\lambda  & \cdots & -w_{2,n-1} & -w_{2,n} \\
			\vdots & \vdots & \ddots & \vdots & \vdots\\
			-w_{1,n-1} & -w_{2,n-1} & \cdots &\alpha_{n-1}-\lambda & -w_{n-1,n}\\
			-\lambda & -\lambda & \cdots & -\lambda &-\lambda
			\end{matrix}\right|,
			\end{equation*}}%
		where the expression in (\ref{eq:alpha}) is applied.

		Using a similar argument, adding the last column to all other columns then leads to (\ref{eq:pencil}). As the eigenvalues of $\mathcal{L}$ are determined by the roots of $\left| \mathcal{L} -  \lambda I_n \right| = 0$, we can assign the spectra of $\mathcal{L}$ by manipulating the weights $w_{i,j}$.  
		
	\begin{figure*}[h]
	{ 
		\begin{equation} \label{eq:pencil}
		\begin{split}
		\left| \mathcal{L}-\lambda I \right| &= 
		\left|  \begin{matrix}
		\alpha_1+w_{1,n}-\lambda & w_{1,n}-w_{1,2} & \cdots & w_{1,n}-w_{1,n-1} & -w_{1,n} \\
		w_{2,n}-w_{1,2} & \alpha_2+w_{2,n}-\lambda  & \cdots & w_{2,n}-w_{2,n-1} & -w_{2,n} \\
		\vdots & \vdots & \ddots & \vdots & \vdots\\
		w_{n-1,n}-w_{1,n-1} & w_{n-1,n}-w_{2,n-1} & \cdots &\alpha_{n-1}+w_{n-1,n}-\lambda & -w_{n-1,n}\\
		0 & 0 & \cdots & 0 &-\lambda
		\end{matrix}\right|.
		\end{split}
		\end{equation}
	 			\hrulefill
	}
\end{figure*}
		
		When $n=2$, we have a special case, and therefore it is considered separately. Equation (\ref{eq:pencil}) becomes 
		{\scriptsize
			\begin{equation*}
			\left|  \mathcal{L} - \lambda I_2 \right| = \left|  \begin{matrix}
			\alpha_1+w_{1,2}-\lambda & -w_{1,2}  \\
			0 &-\lambda
			\end{matrix}\right| = \left|  \begin{matrix}
			2w_{1,2}-\lambda & -w_{1,2}  \\
			0 &-\lambda
			\end{matrix}\right|.
			\end{equation*}}%
		To match the eigenvalues $0$ and $\lambda_1$, we let $w_{1,2}=0.5\lambda_1$, which yields a Laplacian matrix as 
		\begin{equation}
		\mathcal{L} = \begin{bmatrix}
		0.5\lambda_1 & -0.5\lambda_1 \\
		-0.5\lambda_1 & 0.5\lambda_1
		\end{bmatrix},
		\end{equation}
		and proves the desired result for $n=2$.
		
		We continue the proof for the case $n>2$. To match the eigenvalues of $\mathcal{L}$ with the desired ones in (\ref{eq:lamdas}), we let the off-diagonal entries in the lower triangular part of the determinant in (\ref{eq:pencil}) be zero and use the diagonal entries to match the eigenvalues $\lambda_i$ ($i=1,2,\cdots,n$).
		Specifically, the weights $w_{i,j}$ in (\ref{eq:GeneralFormLaplacian}) need to satisfy
		\begin{equation} \label{eqs:eq}
		\left\{
		\begin{split}
		w_{2,n}&=w_{1,2},  \\
		w_{3,n}&=w_{1,3} =w_{2,3},\\
		w_{4,n}&=w_{1,4} = w_{2,4} = w_{3,4},\\
		& \ \ \vdots \\
		w_{n-1,n}&= w_{1,n-1} = w_{2,n-1} = \cdots = w_{n-2,n-1},
		\end{split}\right.
		\end{equation}
		and
		\begin{equation} \label{eqs:sum}
		\alpha_{i} + w_{i,n} = \lambda_{i}, \ \forall  i \in \{1,2,\cdots,n-1\}.
		\end{equation} 
		Hereafter we prove that the equations (\ref{eqs:eq}) and (\ref{eqs:sum}) produce a unique set of nonnegative real weights $w_{i,j}$, which is a necessary and sufficient property to allow for interpretation as a Laplacian matrix, see Remark \ref{rem:StruCond}.

		
		%
		%
		
		For simplicity, we denote
		\begin{equation} \label{eq:aiDef0}
		a_l = w_{n-l,n}, \ l=1,2,\cdots,n-1.
		\end{equation}
		For any $1 \leq l \leq n-2$, it follows from (\ref{eqs:eq}) and the symmetry of $\mathcal{L}$ that 
		\begin{equation} \label{eq:aiDef}
		a_l = w_{k,n-l} = w_{n-l,k}, \ \forall k \in \{1,\cdots,n-l-1\}.
		\end{equation} 
		Furthermore, denote the sum of the above series as
		\begin{equation} \label{eq:SiDef}
		S_l := \sum\limits_{k=1}^{l} a_k, \ l=1,2,\cdots,n-1.
		\end{equation}
		From (\ref{eqs:sum}) and the expression (\ref{eq:alpha}), we have 
		{\begin{align} \label{eq:lambda_i}
			\nonumber \lambda_i & = (w_{i,1}+\cdots+w_{i,i-1}) \\ 
			\nonumber & \ \ \ \ \ +   (w_{i,i+1}+\cdots+w_{i,n-1})+2w_{i,n}
			\\
			\nonumber & = (i-1) a_{n-i} + (a_{n-i-1}+\cdots+a_{1}) + 2a_{n-i}\\
			& = (i+1) a_{n-i} + S_{n-i-1},
			\end{align}}%
		for $i=1,2,\cdots,n-2$. Here, the first equality follows from (\ref{eq:aiDef0}) and (\ref{eq:aiDef}) (with $i = n-l$ for the first term). The latter equation is the result of (\ref{eq:SiDef}).
		
		Rewriting (\ref{eq:lambda_i}) for $l = n-i$ leads to
		\begin{equation} \label{eq:ai}
		a_{l}
		= \dfrac{1}{n-l+1}\left(\lambda_{n-l}-S_{l-1}\right).
		\end{equation}
		
		Now, we prove that $a_l>0$,  $\forall l\in \{1,2,\cdots,n-1\}$. 
		To do so, we consider the cases $l=1$ and $l=2$ explicitly and then proceed by induction. 
		
		For $l=1$, it follows from (\ref{eq:alpha}) and the last equation in (\ref{eqs:eq}) that (\ref{eqs:sum}) can be written as $n w_{n-1,n} = \lambda_{n-1}$, which leads to
		\begin{equation}  
		a_1 = \dfrac{\lambda_{n-1}}{n} = S_1 > 0,
		\end{equation} 
		by the definitions in (\ref{eq:aiDef0}) and (\ref{eq:SiDef}).

		For $l=2$, (\ref{eq:ai}) gives
		\begin{align}\label{eq:solution2}
		\nonumber a_2 & = \dfrac{1}{n-1}\left(\lambda_{n-2}-S_1\right) \\
		& \geq \dfrac{1}{n-1}\left(\lambda_{n-1}-S_1\right) = \dfrac{\lambda_{n-1}}{n} > 0,
		\end{align} 
		where the inequality follows from the ordering of the eigenvalues in (\ref{eq:lamdas}). Then, using (\ref{eq:SiDef}), it follows that
		\begin{equation} 
		S_2 = S_1 + a_2
		= \dfrac{\lambda_{n-2}}{n-1}+\dfrac{(n-2)\lambda_{n-1}}{n(n-1)}.
		\end{equation}
		Note that $\forall m \ne n,  n \ne 0$, we have
		\begin{equation} \label{eq:m+1/n}
		\dfrac{1}{n-m}+\dfrac{m(n-m-1)}{n(n-m)} = \dfrac{m+1}{n}.
		\end{equation}
		Using the above equation with $m=1$ and the inequality $\lambda_{n-2} \geq \lambda_{n-1}$, we show bounds on $S_2$ as
		\begin{equation} \label{eq:S2}
		\begin{split}
		S_2 & \geq \left[\dfrac{1}{n-1}+\dfrac{(n-2)}{n(n-1)}\right] \lambda_{n-1} = \dfrac{2\lambda_{n-1}}{n}, \\
		S_2 &\leq \left[\dfrac{1}{n-1}+\dfrac{(n-2)}{n(n-1)}\right] \lambda_{n-2} = \dfrac{2\lambda_{n-2}}{n}.
		\end{split}
		\end{equation}
		To proceed with induction on $l$ for $l>2$, we assume both $a_l > 0$ and
		\begin{equation} \label{eq:Si}
		\dfrac{l\lambda_{n-1}}{n} \leq S_l \leq \dfrac{l\lambda_{n-l}}{n},
		\end{equation}
		for  $2 < l < n-1$. 
		Then, we obtain from (\ref{eq:ai}) and (\ref{eq:Si}) that 
		\begin{align} \label{eq:a_i+1}
		\nonumber a_{l+1}
		&= \dfrac{1}{n-l}\left(\lambda_{n-l-1}-S_{l}\right) \\ &\geq  \dfrac{1}{n-l}\left(\lambda_{n-l-1}-\dfrac{l\lambda_{n-l}}{n}\right)
		\geq
		\dfrac{\lambda_{n-l}}{n}>0,
		\end{align}
		after which the first line in (\ref{eq:a_i+1}) yields
		\begin{equation}
		S_{l+1} = S_l + a_{l+1}
		= \dfrac{\lambda_{n-l-1}}{n-l}+\dfrac{(n-l-1)S_l}{n-l}.
		\end{equation}
		The upper and lower bounds on $S_{l+1}$ are implied by (\ref{eq:Si}) as
		\begin{equation}
		\begin{split}
		S_{l+1} &\geq \dfrac{\lambda_{n-l-1}}{n-l}+\dfrac{l(n-l-1)\lambda_{n-1}}{(n-l)n}, \\  S_{l+1} &\leq \dfrac{\lambda_{n-l-1}}{n-l}+\dfrac{l(n-l-1)\lambda_{n-l}}{(n-l)n}.
		\end{split}
		\end{equation}
		Using the relation $\lambda_{n-l-1} \geq \lambda_{n-l} \geq \lambda_{n-1}$ and the equation (\ref{eq:m+1/n}) with $m=l$, we obtain
		\begin{equation}
		\dfrac{(l+1)\lambda_{n-1}}{n} \leq S_{l+1} \leq \dfrac{(l+1)\lambda_{n-l-1}}{n}.
		\end{equation}
		Consequently, by induction, we now verify that $a_l>0$, $\forall l \in \{1,2,\cdots,n-1\}$. As the parameters $a_l$ uniquely characterize all the the weights $w_{i,j}$ in (\ref{eq:GeneralFormLaplacian}) through (\ref{eq:aiDef0}) and (\ref{eq:aiDef}), it follows that $w_{i,j} > 0 $ for all $(i,j)$.
		
		In summary, there always exist a set of weights $w_{i,j} > 0$ such that $\mathcal{L}$ in (\ref{eq:GeneralFormLaplacian}) has the eigenvalues matching the desired spectrum $\lambda_1  \geq \cdots \geq \lambda_{n-1} > \lambda_{n}=0$. The matrix $\mathcal{L}$ satisfies all properties stated in Remark \ref{rem:StruCond} and thus is a Laplacian matrix representing an undirected graph.
		Therefore, we conclude that if $\mathcal{N}$ is diagonalizable and has a single zero eigenvalue while all the other eigenvalues are real positive, then there always exists a similarity transformation between $\mathcal{N}$ and a Laplacian matrix. This finalizes the proof of Theorem~\ref{thm:LapReal}.
\end{pf}


\bibliographystyle{abbrv}        
\bibliography{NetworkReduction}

\end{document}